\documentclass[reqno,a4paper]{amsart}

\usepackage{amsmath,amsfonts,amssymb}
\usepackage{verbatim}
\usepackage{enumerate}
\usepackage{url}
\usepackage{tikz}
\usepackage{setspace}
\usepackage[left=2.5cm,right=2.5cm, top=3cm, bottom=2.5cm]{geometry}
\usepackage[utf8]{inputenc} 
\usepackage[T1]{fontenc}
\usepackage{hyperref}
\usetikzlibrary{arrows}

\def\N{\mathbb N}

\def\Z{\mathbb Z}

\theoremstyle{plain}
\newtheorem{theorem}{Theorem}[section]
\newtheorem{lemma}[theorem]{Lemma}

\newtheorem{corollary}[theorem]{Corollary}

\newtheorem{remark}[theorem]{Remark}

\def\proof{\noindent {\it Proof: }}
\def\qed{\hfill\hbox{$\square$}}

\theoremstyle{definition}

\numberwithin{equation}{section}

\subjclass[2010]{11B50 (primary), 11P70 (secondary)}
\title{The subsums of zero-sum free sequences in finite cyclic groups}
\keywords{Zero-sum problem, Davenport constant, inverse zero-sum problem}

\author[S\'avio Ribas]{S\'avio Ribas}
\address{
Instituto Federal de Educa\c{c}\~ao, Ci\^encia e Tecnologia de Minas Gerais\\
IFMG - Campus Ouro Preto\\
Ouro Preto, MG\\
35400-000\\
Brazil\\
}
\email{savio.ribas@ifmg.edu.br }

\date{\today}

\onehalfspace

\begin{document}

\maketitle

\begin{abstract}
Let $\Z_n$ be the cyclic group of order $n \ge 3$ additively written. S. Savchev \& F. Chen (2007) proved that for each zero-sum free sequence $S = a_1 \bullet \dots \bullet a_t$ over $\Z_n$ of length $t > n/2$, there is an integer $g$ coprime to $n$ such that, if $\overline{r}$ denotes the least positive integer in the congruence class $r$ modulo $n$, then $\sum_{i=1}^t \overline{ga_i} < n$. Under the same hypothesis, in this paper we show that 
$$\left\{ \sum_{i \in \Lambda} \overline{ga_i} \;\; \Bigg| \;\; \Lambda \subset \{1,2,\dots,t\} \right\} = \left\{ 1, 2, \dots, \sum_{i=1}^t \overline{ga_i}\right\}.$$ It simplifies many calculations on inverse zero-sum problems.
\end{abstract}

\section{Introduction}

Given a finite group $G$, the {\em Zero-Sum Problems} study conditions to ensure that a given sequence over $G$ has a non-empty subsequence with prescribed properties (such as length, repetitions, weights) such that the product of its elements is equal to the identity of $G$. This class of problems have been extensively studied for abelian groups (see, for example, the surveys \cite{Car, GaGe}), and since little more than a decade there are some results for non-abelian groups (see, for example, \cite{Bas, MR1, MR2, GeGr, Gr, Oh, Oh2, OlWh, ZhGa}).

Denote by $[a,b]$ the interval $\{n \in \N; \; a \le n \le b\}$. By a sequence $S$ over $G$ we mean an element $S$ (finite and unordered) of the free abelian monoid $\mathcal F(G)$, equipped with the sequence concatenation product denoted by $\bullet$. A sequence $S \in \mathcal F(G)$ has the form $$S = g_1 \bullet \dots \bullet g_k = \prod_{i=1}^k g_i \in \mathcal F(G)$$ where $g_i \in G$ are the {\em terms} of $S$ and $k = |S| \ge 0$ is the {\em length} of $S$. Since the sequences are unordered, $$S = \prod_{i=1}^k g_i = \prod_{i=1}^k g_{\tau(i)} \in \mathcal F(G)$$ for any permutation $\tau: [1,k] \to [1,k]$. Given $g \in G$ and $t \ge 0$, we abbreviate $$g^{[t]} = \underbrace{g \bullet \dots \bullet g}_{t \text{ times}}.$$
For $g \in G$, the {\em multiplicity} of the term $g$ in $S$ is denoted by $v_g(S) = \#\{i \in [1,k]; \; g_i = g\}$, therefore our sequence $S$ may also be written as $$S = \prod_{g \in G} \; g^{[v_g(S)]}.$$
A sequence $T$ is called a {\em subsequence} of $S$ if $T | S$ in $\mathcal F(G)$, or equivalently, $v_g(T) \le v_g(S)$ for all $g \in G$. From now on, we assume that $G$ is abelian and additively written. \\
We also define:
\begin{align*}
\sigma(S) &= g_1 + \dots + g_k \in G \;\;\; \text{ the {\em sum} of $S$}; \\
\Sigma(S) &= \bigcup_{T | S \atop |T| \ge 1} \{\sigma(T)\} \subset G \;\;\; \text{ the {\em set of subsequence sums} of $S$}; \\
S \bullet T^{-1} &= \prod_{g \in G \atop {v_g(S) > 0 \atop v_g(T) = 0}} g \;\;\; \text{ the subsequence of $S$ formed by the terms that do not lie in $T$}; \\
S \cap K &= \prod_{g \in S \atop g \in K} g \;\;\; \text{ the subsequence of $S$ that lie in a subset $K$ of $G$}. 
\end{align*}
The sequence $S$ is called
\begin{itemize}
\item[-] {\em zero-sum free} if $0 \not\in \Sigma(S)$;
\item[-] {\em zero-sum sequence} if $\sigma(S) = 0$.
\end{itemize}

An important type of zero-sum problem is to determine the {\em Davenport constant} of a finite group $G$: This constant, denoted by $d(G)$, is the maximal integer such that there exists a sequence over $G$ (repetition allowed) of length $d(G)$ which is zero-sum free, ie,
$$d(G) = \sup\{|S|>0; \; S \in \mathcal F(G) \text{ is zero-sum free}\}.$$ 

Using Pigeonhole Principle on the partial sums $g_1, g_1+g_2, \dots, g_1+\dots+g_n$, it is easy to show that if $S \in \mathcal F(\Z_n)$ and $|S| = n$ then $S$ is not zero-sum free. Furthermore, if $\gcd(a,n)=1$ then the sequence $a^{[n-1]}$ is zero-sum free. Thus $d(\Z_n) = n-1$. The zero-sum free sequence $S = (1,0,\dots,0)^{[n_1-1]} \bullet \dots \bullet (0,\dots,0,1)^{[n_k-1]}$ over $\Z_{n_1} \oplus \dots \oplus \Z_{n_k}$ shows that 
\begin{equation}\label{lowerboundabelian}
d(\Z_{n_1} \oplus \dots \oplus \Z_{n_k}) \ge \sum_{i=1}^k (n_i - 1).
\end{equation}
Olson \cite{Ols1, Ols2} proved that the equality holds for $p$-groups and in the case of rank two. It is conjectured that the equality also holds for groups of the form $\Z_n^k$. However, the equality does not hold for a general abelian group of order $k \ge 4$ (see, for example, \cite{GeSc, Maz}).

Denote by $\overline{a}$ the least positive integer in the congruence class $a \pmod n$ and let $n \ge 2k + 1 \ge 3$. Bovey, Erd\H{o}s \& Niven \cite{BEN} proved that if $S \in \mathcal F(\Z_n)$ is zero-sum free and $|S| = n - k$, then $$v_a(S) \ge n - 2k + 1$$ for some $a | S$. Moreover, they showed that this inequality is the best possible whether $n \ge 3k - 2$. Savchev \& Chen \cite{SC} proved that if $2k + 1 \le n \le 3k - 2$ then there exists $a | S$ such that
$$v_a(S) \ge n - k - \left\lfloor \frac{n-1}{3} \right\rfloor = 
\begin{cases}
\frac{2n+3}3 - k \quad \text{ if $n \equiv 0 \pmod 3$;} \\
\frac{2n+1}3 - k \quad \text{ if $n \equiv 1 \pmod 3$;} \\
\frac{2n+2}3 - k \quad \text{ if $n \equiv 2 \pmod 3$.}
\end{cases}$$
In addition, they proved that for each zero-sum free sequence $S = a_1 \bullet \dots \bullet a_{n-k} \in \mathcal F(\Z_n)$ there is an integer $g$ with $\gcd(g,n)=1$ such that $$\sum_{i=1}^{n-k} \overline{ga_i} \le n - 1.$$

In this paper, we show that if $n \ge 2k+1 \ge 3$ then 
$$\left\{ \; \sum_{i \in \Lambda} \overline{ga_i} \;\; \Bigg| \;\; \Lambda \subset [1,t] \; \right\} = \left[ \; 1 \; , \; \sum_{i=1}^t \overline{ga_i} \; \right].$$ 

The proof is done in Section \ref{proof}. In Section \ref{lemas} we present some auxiliary lemmas, while in Section \ref{application} we present an application that simplifies some calculations on inverse problems.

\vspace{2mm}

\section{Preliminary}\label{lemas}

Our proof is based on the following results:

\begin{lemma}\label{lema1}
Let $S \in \mathcal F(\Z)$ of the form $S = 1^{[v_1]} \bullet 2^{[v_2]}$, where $v_1, v_2 \ge 1$. Then $\Sigma(S) = [ 1 , v_1 + 2v_2 ]$.
\end{lemma}

\proof
The inclusion $\subset$ is clear. For the opposite $\supset$, let $N \in [1,v_1+2v_2]$. Then $N = \alpha + 2\beta$, where $0 \le \alpha \le v_1$ and $0 \le \beta \le v_2$ (but not $\alpha = \beta = 0$), which implies $N = \sum_{i=1}^{\alpha} 1 + \sum_{i=1}^{\beta} 2 \in \Sigma(S)$.
\qed

\begin{lemma}\label{lema2}
Let $S \in \mathcal F(\Z)$ of the form $S = 1^{[v_1]} \bullet 3^{[v_3]}$, where $v_1 \ge 2$ and $v_3 \ge 1$. Then $\Sigma(S) = [ 1 , v_1 + 3v_2 ]$.
\end{lemma}

\proof
The inclusion $\subset$ is clear. For the opposite $\supset$, let $N \in [1,v_1+3v_2]$. Then $N = \alpha + 3\gamma$, where $0 \le \alpha \le v_1$ and $0 \le \gamma \le v_2$ (but not $\alpha = \gamma = 0$), which implies $N = \sum_{i=1}^{\alpha} 1 + \sum_{i=1}^{\gamma} 3 \in \Sigma(S)$.
\qed

\begin{lemma}\label{lema3}
Let $G$ be a finite abelian group and let $S \in \mathcal F(G)$. If $T|S$ then $$\sigma(T) + \sigma(S \bullet T^{-1}) = \sigma(S).$$
In particular, if $\sigma(T) = t$ then $\sigma(S \bullet T^{-1}) = \sigma(S) - t$.
\end{lemma}

\proof
Trivial.
\qed

\vspace{2mm}

\section{The subsums of a zero-sum free sequence}\label{proof}

\begin{theorem}
Let $n,k$ be positive integers with $n \ge 2k + 1 \ge 3$ and let $S = a_1 \bullet \dots \bullet a_{n-k} \in \mathcal F(\Z_n)$ be a zero-sum free sequence with 
\begin{equation}\label{savche}
\sum_{i=1}^{n-k} \overline{a_i} \le n-1.
\end{equation}
Then $$\Sigma(S) = [1,\sigma(S)].$$
In other words, for all $1 \le t \le \displaystyle \sum_{i=1}^{n-k} \overline{a_i}$ there exist $a_{i_1} \bullet \dots \bullet a_{i_r} \; | \; S$ such that $\displaystyle \sum_{j=1}^r \overline{a_{i_j}} = t$.
\end{theorem}

\proof
The inclusion $\subset$ is clear. For the opposite $\supset$, by Lemma \ref{lema3} and inequality \ref{savche}, it is enough show that 
\begin{equation}\label{n/2}
\left[ 1, \left\lfloor \frac{n}{2} \right\rfloor \right] \subset \Sigma(S).
\end{equation}

For simplicity, we denote $v_i(S)$ by $v_i$. Suppose that $v_1 \le 1$. Then $$\sum_{i=1}^{n-k} \overline{a_i} \ge 1 + 2(n-k-1) = 2n - 2k - 1 \ge n,$$
since $n \ge 2k + 1$, which is an absurd. Therefore, $v_1 \ge 2$.

Let $v_{\ell} = \max\{v_i\}$ (that is, $\ell$ is the term of greatest multiplicity in $S$). By the results of Bovey, Erd\H{o}s \& Niven \cite{BEN} and Savchev \& Chen \cite{SC}, we split into two cases:

\begin{enumerate}[(i)]
\item {\bf Case $n \ge 3k - 2$:} we have $v_{\ell} \ge n - 2k + 1 \ge \left\lceil \frac{n-1}{3} \right\rceil$.

\begin{enumerate}[{(i.}1)]
\item If $\ell \ge 3$ then 
$$\sum_{i=1}^{n-k} \overline{a_i} > 3v_{\ell} \ge 3 \left( \frac{n-1}{3} \right) = n-1, \quad \text{contradiction!}$$

\item If $\ell = 2$, since 
$$\sum_{i=1}^{n-k} \overline{a_i} > 2v_2 \ge 2\left\lceil \frac{n-1}{3} \right\rceil \ge \frac n2$$ 
and $v_1 \ge 2$, Lemma \ref{lema1} implies that $\left[1, \left\lfloor \frac{n}{2} \right\rfloor \right] \subset \Sigma(S)$. Therefore, inclusion (\ref{n/2}) holds and we are done in this subcase.

\item If $\ell = 1$, define $0 \le \alpha \le k-1$ by equation $v_1 = n - 2k + 1 + \alpha$. By Pigeonhole Principle, the average of the elements $\overline{a_i}$ without the terms $1$ and $2$ must be at least $3$ (otherwise, there would be $\overline{a_j} < 3$ that could be removed before). Hence,
\begin{equation}\label{pcp}
3 \le \frac{\left( \displaystyle\sum_{i=1 \atop \overline{a_i} \ge 3}^{n-k} \overline{a_1} \right)}{n - k - v_1 - v_2} \le \frac{n - 1 - v_1 - 2v_2}{n - k - v_1 - v_2}
\end{equation}
which implies 
$$v_2 \ge
\begin{cases}
k - 1 - 2\alpha \quad \text{ if $\alpha \le \frac{k-1}{2}$} \\
\quad \quad 0 \quad \quad \quad \, \text{ if $\alpha > \frac{k-1}{2}$} 
\end{cases}.$$

\noindent If $\alpha \le \frac{k-1}{2}$ then we have 
$$\sum_{i=1}^{n-k} \overline{a_i} \ge v_1 + 2v_2 \ge n-1 - 3 \left( \frac{k-1}{2} \right) \ge \frac{n}{2},$$
and we are done by Lemma \ref{lema1} and inclusion (\ref{n/2}).

\noindent Otherwise, if $\alpha > \frac{k-1}{2}$ then we have 
$$v_1 \ge n - 2k + 1 + \frac{k}{2} \ge \frac{n}{2},$$
thus every element in $\left[ 1, \left\lfloor \frac{n}{2} \right\rfloor \right]$ can be obtained as sum of $1$'s and we are done in this case by inclusion (\ref{n/2}).
\end{enumerate}

\vspace{2mm}

\item {\bf Case $2k + 1 \le n < 3k - 2$:} we have $v_{\ell} \ge n - k - \lfloor \frac{n-1}{3} \rfloor \ge \frac{n+5}{6}$.

\begin{enumerate}[{(ii.}1)]
\item If $\ell \ge 4$ then 
$$\sum_{i=1}^{n-k} \overline{a_i} \ge 4v_{\ell} + \sum_{i=1 \atop i \neq \ell}^{n-k} \overline{a_i} \ge 4 \left( \frac{n+5}{6} \right) + \left\lfloor \frac{n-1}{3} \right\rfloor > n-1, \quad \text{contradiction!}$$

\item If $\ell = 3$, since 
$$\sum_{i=1}^{n-k} \overline{a_i} > 3v_3 \ge \frac{n+5}{2} > \frac n2$$ 
and $v_1 \ge 2$, Lemma \ref{lema2} implies that $\left[1, \left\lfloor \frac{n}{2} \right\rfloor \right] \subset \Sigma(S)$. Therefore, inclusion (\ref{n/2}) holds and we are done in this subcase.

\item If $\ell = 2$, since $v_1 \ge 2$ and $$\sum_{i=1}^{n-k} \overline{a_i} \ge 2v_2 + \sum_{i=1 \atop i \neq 2}^{n-k} \overline{a_i} \ge 2 \left( \frac{n+5}{6} \right) + \left\lfloor \frac{n-1}{3} \right\rfloor > \frac{n-1}{2},$$
then Lemma \ref{lema1} implies that $\left[1, \left\lfloor \frac{n}{2} \right\rfloor \right] \subset \Sigma(S)$. Therefore, inclusion (\ref{n/2}) holds and we are done in this subcase.

\item If $\ell = 1$, define $0 \le \alpha \le \left\lfloor \frac{n-1}{3} \right\rfloor$ by equation $v_1 = n - k - \left\lfloor \frac{n-1}{3} \right\rfloor + \alpha$. Using Pigeonhole Principle again, the average of the elements $\overline{a_i}$ without the terms $1$ and $2$ must be at least $3$, hence inequality (\ref{pcp}) implies that 
$$v_2 \ge
\begin{cases}
2 \left\lfloor \frac{n-1}{3} \right\rfloor - 2\alpha + k + 1 \quad \text{ if $\alpha \le \left\lfloor \frac{n-1}{3} \right\rfloor + \frac{k+1}{2}$} \\
\quad \quad \quad \quad \quad 0 \quad \quad \quad \quad \;\; \, \text{ if $\alpha > \left\lfloor \frac{n-1}{3} \right\rfloor + \frac{k+1}{2}$} 
\end{cases}.$$

\noindent If $\alpha \le \left\lfloor \frac{n-1}{3} \right\rfloor + \frac{k+1}{2}$ then we have 
$$\sum_{i=1}^{n-k} \overline{a_i} \ge v_1 + 2v_2 \ge n - \frac{k}{2} + 1 \ge \frac{n}{2},$$
and we are done by Lemma \ref{lema1} and inclusion (\ref{n/2}).

\noindent Otherwise, if $\alpha > \frac{k-1}{2}$ then we have 
$$v_1 \ge n - \frac{k - 1}{2} \ge \frac{n}{2},$$
thus every element in $\left[ 1, \left\lfloor \frac{n}{2} \right\rfloor \right]$ can be obtained as sum of $1$'s and we are done in this case by inclusion (\ref{n/2}).
\end{enumerate}
\end{enumerate}
\qed

{\remark By the results of Savchev \& Chen, there is no loss of generality in the assumption (\ref{savche}).}

{\remark The hypothesis $n \ge 2k+1 \ge 3$ can not be removed. In fact, if $n=5$, $k=3$, $S_1 = 1 \bullet 3 \in \mathcal F(\Z_5)$ and $S_2 = 2^{[2]} \in \mathcal F(\Z_5)$, then $2 \not\in \Sigma(S_1) = \{1,3,4\}$ and $1,3 \not\in \Sigma(S_2) = \{2,4\}$.}

\vspace{2mm}

\begin{corollary}\label{corol}
Let $n,k$ be positive integers with $n \ge 2k + 1 \ge 3$ and let $S = a_1 \bullet \dots \bullet a_{n-k} \in \mathcal F(\Z_n)$ be a zero-sum free sequence satisfying Equation \ref{savche}. Then $$[1,n-k] \subset \Sigma(S).$$
In other words, for all $1 \le t \le n-k$ there exist $a_{i_1} \bullet \dots \bullet a_{i_r} \; | \; S$ such that $\displaystyle \sum_{j=1}^r \overline{a_{i_j}} = t$.
\end{corollary}

\section{Simplifying calculations}\label{application}

In this section, we present an example of application that simplify some calculations on inverse zero-sum problems. Let $G$ be a finite non-abelian group multiplicatively written, and $S = g_1 \bullet \dots \bullet g_k \in \mathcal F(G)$. Define:
\begin{align*}
\pi(S) &= \{ g_{\tau(1)} \cdot \;\! \dots \;\! \cdot g_{\tau(k)} \in D_{2n} ; \; \tau \text{ is a permutation in } [1,k] \} \;\;\; \text{ the {\em set of products} of $S$}; \\
\Pi(S) &= \bigcup_{T | S \atop |T| \ge 1} \pi(T) \subset D_{2n} \;\;\; \text{ the {\em set of subsequence products} of $S$}.
\end{align*}
The sequence $S$ is called {\em product-one free} if $1 \not\in \Pi(S)$. The {\em small Davenport constant} $d(G)$ can be defined analogously for $G$:
$$d(G) = \sup\{|S|>0; \; S \in \mathcal F(G) \text{ is product-one free}\}.$$ 

Let $D_{2n}$ be the Dihedral Group, ie, the group generated by $x$ and $y$ satisfying $x^2 = y^n = 1$ and $yx = xy^{-1}$. Zhuang \& Gao \cite{ZhGa} proved $d(D_{2n}) = n$. In a joint work with Brochero Martínez \cite{MR2}, we exhibit all the extremal length product-one free sequences over $D_{2n}$, showing that:

\begin{theorem}[\cite{MR2}, Theorem 1.3]\label{inversedihedral}
Let $n \ge 3$ and $S \in \mathcal F(G)$ such that $|S| = n$.
\begin{enumerate}
\item If $n \ge 4$ then $S$ is product-one free if and only if for some $1 \le t \le n-1$ with $\gcd(t,n)=1$ and $0 \le s \le n-1$, $S = (y^t)^{[n-1]} \bullet xy^s$.
\item If $n=3$ then $S$ product-one free if and only if either $S = x \bullet xy \bullet xy^2$ or $S = (y^t)^{[2]} \bullet xy^{\nu}$ for $t \in \{1,2\}$ and $\nu \in \{0,1,2\}$.
\end{enumerate}
\end{theorem}

In the original paper, the proof was obtained considering the cases 
\begin{itemize}
\item $|S \cap \langle y \rangle| = n$,
\item $|S \cap \langle y \rangle| = n-1$,
\item $|S \cap \langle y \rangle| = n-2$,
\item $|S \cap \langle y \rangle| = n-3$,
\item $n - 2 \lfloor \log_2 n \rfloor - 1 \le |S \cap \langle y \rangle| \le n-4$, and
\item $|S \cap \langle y \rangle| \le n - 2 \lfloor \log_2 n \rfloor - 2$,
\end{itemize}
besides considering the initial cases $3 \le n \le 7$ separately. Only the third and fourth cases resemble each other. The case $n$ prime had already been done in \cite{MR1}, which reduced the handwork at least in the initial steps.

\vspace{2mm}

\proof
In our new approach, there is no need to consider any initial steps, and in addition, the first two cases keep identical to the original paper. We propose to consider only the following 
\begin{itemize}
\item $|S \cap \langle y \rangle| = n$,
\item $|S \cap \langle y \rangle| = n-1$,
\item $\frac{n}{2} < |S \cap \langle y \rangle| \le n-2$, and 
\item $|S \cap \langle y \rangle| \le \frac{n}{2}$.
\end{itemize}

For the third case, we notice that $xy^{\alpha} \bullet xy^{\beta} | S$ for some $\alpha \not\equiv \beta \pmod n$ (otherwise, $xy^{\alpha} \cdot xy^{\beta} = 1$). Assume that $\overline{\alpha - \beta} \in \left[ 1, \left\lfloor \frac{n}{2} \right\rfloor \right]$ (otherwise, switch $\alpha$ and $\beta$). Using the result of Savchev \& Chen if needed, we may assume that if $S \cap \langle y \rangle$ does not satisfy the hypotheses of Corollary \ref{corol} then $S$ is not product-one free. This corollary ensures that $\overline{\alpha - \beta}$ can be obtained as sum of the exponents of the elements in $S \cap \langle y \rangle$, say, $\overline{\alpha - \beta} = {\gamma_1} + {\gamma_2} + \dots + {\gamma_r}$, where $y^{\gamma_i} | S$. Therefore, $xy^{\alpha} \cdot xy^{\beta} \cdot y^{\gamma_1} \cdot y^{\gamma_2} \dots y^{\gamma_r} = 1$ and $S$ is not product-one free.

For the last case, define $k \in \N$ by the equation $|S \cap \langle y \rangle| = n - k$, where $k \ge \frac{n}{2}$. Let $T = S \bullet (S \cap \langle y \rangle)^{-1}$. Then $|T| = k$, which implies that $T$ produces $\left\lfloor \frac k2 \right\rfloor$ new elements in $S \cap \langle y \rangle$. Let $xy^{\alpha} \cdot xy^{\beta} y = y^{\beta - \alpha}$ one of these, where $1 \le \alpha - \beta \le \frac{n}{2}$. If $n \ge k+1$ then $n - k + \left\lfloor \frac k2 \right\rfloor - 1 \ge \left\lfloor \frac n2 \right\rfloor$, thus Corollary $\ref{corol}$ ensures that the $n-k$ terms of $S \cap \langle y \rangle$ joint with the $\left\lfloor \frac{k}{2} \right\rfloor - 1$ new elements from $T$ yields a product of the form $P = y^{\alpha - \beta}$. Therefore, $P \cdot xy^{\alpha} \cdot xy^{\beta} = 1$, and $S$ is not product-one free. It only remains the case $n = k$. But in this case the elements must be all distinct (otherwise the product of a identical pair would be $1$), therefore $S = x \bullet xy \bullet xy^2 \bullet \dots \bullet xy^{n-1}$. If $n \ge 4$ then $x \cdot xy \cdot xy^3 \cdot xy^2 = 1$ and $S$ is not product-one free. If $n=3$, the sequence $x \bullet xy \bullet xy^2$ over $D_6$ is product-one free, and we are done.
\qed


\vspace{2mm}



\begin{thebibliography}{99}

\bibitem{Bas} J. Bass;
{\em Improving the Erd\H{o}s-Ginzburg-Ziv theorem for some non-abelian groups.}
J. Number Theory {\bf 126} (2007), 217-236.

\bibitem{BEN} J.D. Bovey, P. Erd\H{o}s, I. Niven;
{\em Conditions for a zero sum modulo n.}
Canad. Math. Bull. Vol. {\bf 18} (1), (1975), 27-29.

\bibitem{MR1} F.E. Brochero Mart\'inez, S\'avio Ribas;
{\em Extremal product-one free sequences in $C_q \rtimes_s C_m$.}
Submitted for publication (2017). Available at \url{https://arxiv.org/pdf/1610.09870.pdf}.

\bibitem{MR2} F.E. Brochero Mart\'inez, S\'avio Ribas;
{\em Extremal product-one free sequences in Dihedral and Dicyclic Groups.}
Discrete Mathematics {\bf 341} (2018), 570-578.

\bibitem{Car} Y. Caro;
{\em Zero-sum problems - A survey.}
Discrete Mathematics {\bf 152} (1996) 93-113.

\bibitem{GaGe2} W. Gao, A. Geroldinger;
{\em On long minimal zero sequences in finite abelian groups.}
Period. Math. Hungar. {\bf 38 (3)} (1999), 179-211. 

\bibitem{GaGe3} W. Gao, A. Geroldinger;
{\em On zero-sum sequences in $\Z/n\Z \oplus \Z/n\Z$.}
Integers: Electronic J. of Comb. Number Theory {\bf 3} (2003), \#A8.

\bibitem{GaGe} W. Gao, A. Geroldinger;
{\em Zero-sum problems in finite abelian groups: a survey.}
Expo. Math. {\bf 24} (2006), 337-369.

\bibitem{GaGeGr} W. Gao, A. Geroldinger, D.J. Grynkiewicz;
{\em Inverse zero-sum problems III.}
Acta Arithmetica {\bf 141.2} (2010), 103-152.

\bibitem{GaGeSc} W. Gao, A. Geroldinger, W.A. Schmid;
{\em Inverse zero-sum problems.}
Acta Arithmetica {\bf 128.3} (2007), 245-279.

\bibitem{GeGr} A. Geroldinger, D.J. Grynkiewicz;
{\em The large Davenport constant I: Groups with a cyclic, index 2 subgroup.}
J. of Pure and Applied Algebra {\bf 217} (2013), 863-885.

\bibitem{GeSc} A. Geroldinger, R. Schneider;
{\em On Davenport's constant.}
J. Combin. Theory Ser. A {\bf 61} (1992), 147-152.

\bibitem{Gr} D.J. Grynkiewicz;
{\em The large Davenport constant II: General upper bounds.}
J. of Pure and Applied Algebra {\bf 217} (2013), 2221-2246.

\bibitem{Maz} M. Mazur;
{\em A Note on the Growth of Davenport's constant.}
Manuscripta mathematica {\bf 74.3} (1992), 229-236.

\bibitem{Oh} J.S. Oh;
{\em On the algebraic and arithmetic structure of the monoid of product-one sequences.}
J. Commut. Algebra, to appear.

\bibitem{Oh2} J.S. Oh;
{\em On the algebraic and arithmetic structure of the monoid of product-one sequences II.}
Available at: \url{https://arxiv.org/abs/1802.02851}.

\bibitem{Ols1} J.E. Olson;
{\em A combinatorial problem on finite Abelian groups I.}
J. Number Theory {\bf 1} (1969), 8-10.

\bibitem{Ols2} J.E. Olson;
{\em A combinatorial problem on finite Abelian groups II.}
J. Number Theory {\bf 1} (1969), 195-199.

\bibitem{OlWh} J. Olson, E.T. White; 
{\em Sums from a sequences of group elements.}
Number Theory and Algebra, Academic Press, New York (1977), 215-222.

\bibitem{SC} S. Savchev, F. Chen;
{\em Long zero-free sequences in finite cyclic groups.}
Discrete Mathematics {\bf 307} (2007), 2671-2679.

\bibitem{Sc} W.A. Schmid;
{\em Inverse zero-sum problems II.}
Acta Arithmetica {\bf 143}, no. 4 (2010), 333-343.

\bibitem{ZhGa} J.J. Zhuang, W. Gao;
{\em Erd\H{o}s-Ginzburg-Ziv theorem for dihedral groups of large prime index.}
European J. Combin. {\bf 26} (2005), 1053-1059.

\end{thebibliography}
\end{document}